\documentclass{amsart}

\newtheorem{theorem}{Theorem}[section]

\theoremstyle{definition}

\newtheorem{example}[theorem]{Example}
\newtheorem{corollary}[theorem]{Corollary}
\theoremstyle{remark}
\newtheorem{remark}[theorem]{Remark}

\numberwithin{equation}{section}

\newcommand{\Z}{{\mathbb Z}}
\newcommand{\R}{{\mathbb R}}\newcommand{\T}{{\mathbb T}}
\newcommand{\lau}{{\operatorname{lau}}}\newcommand{\fun}{{\operatorname{fun}}}
\newcommand{\const}{{\operatorname{const}}}
\newcommand{\Del}{{\Delta}}

\begin{document}

\title{A tropical Nullstellensatz}
\author{Eugenii Shustin}
\address{School of Math. Sciences, Tel Aviv University, Ramat Aviv, 69978 Tel Aviv, Israel}
\email{shustin@post.tau.ac.il}
\thanks{The first author was supported by the grant
from the High Council for Scientific and Technological Cooperation
between France and Israel and by the grant no. 465/04 from the
Israel Science Foundation.}
\author{Zur Izhakian}
\address{School of Computer Sciences, Tel Aviv University, Ramat Aviv, 69978 Tel Aviv, Israel}
\email{zzur@post.tau.ac.il}
\thanks{The second author was supported by the grant
from the High Council for Scientific and Technological Cooperation
between France and Israel.} 

\date{August 22, 2005.}


\begin{abstract}
We suggest a version of Nullstellensatz over the tropical
semiring, the real numbers equipped with operations of maximum and
summation.
\end{abstract}

\maketitle

\section{Introduction}\label{intro}

The tropical mathematics is a mathematics over {\it the tropical
semiring}, the real numbers equipped with the operations of
maximum and summation, corresponding to the classical operations
of addition and multiplication, respectively (see in \cite{KM,Vi}
the representation of the tropical semiring as a limit of the
semiring of non-negative real numbers $(\R_+,+,\cdot)$). Sometimes
the tropical semiring is extended by $-\infty$, the neutral
element for the maximum operation, but we shall not need it. We
write the tropical operations in quotes, i.e.,
$$"a+b\ "=\max\{a,b\},\quad "ab\ "=a+b\ .$$

A rapid development of the tropical mathematics over the last
years, especially, of the tropical algebraic geometry, has lead to
spectacular applications in the classical algebraic geometry (see,
for example,
\cite{EKL,I,IKS,IKS2,KS,KT,Mi0,Mi1,Mi,Sh0,SS1,SS,Vi}). The
tropical objects reveal unexpectedly much similarity with
classical objects, which is based on the theory of large complex
limits, logarithmic limits, non-Archimedean valuations, toric
geometry etc. Here we suggest a tropical analogue of
Nullstellensatz (in a bit different context this problem was
stated in \cite{A}, Question 16).

Similarly to the classical case, we understand Nullstellensatz as
a criterion for a polynomial to belong to the radical of an ideal.
Let us give necessary definitions (cf.
\cite{I,Mi1,Mi,Sh0,SS1,Vi}). {\it A tropical polynomial} in $n$
variables is a function $F:\R^n\to\R$ given by
\begin{equation}F(x)="\sum_{\omega\in\Omega}c_\omega x^\omega\
"=\max_{\omega\in \Omega}(\langle x,\omega\rangle+c_\omega)\
,\label{e5}\end{equation} where $\Omega$ is a non-empty finite set
of points in $\Z^n$ with non-negative coordinates,
$\langle*,*\rangle$ denotes the scalar product, and
$c_\omega\in\R$, $\omega\in\Omega$. This is a convex piece-wise
linear function. Its Newton polytope $\Del(F)$ is the convex hull
of the set $\Omega$.

The space of tropical polynomials in $n$ variables is denoted by
$"\R[x_1,...,x_n]\ "$\footnote{The quotation marks mean that we
supply this space by the tropical operations, maximum and
summation.}. For $F_1,...,F_k\in"\R[x_1,...,x_n]\ "$, we define
{\it a tropical polynomial ideal $I(F_1,...,F_k)$, generated by}
$F_1,...,F_k$ in $"\R[x_1,...,x_n]\ "$, as the set of tropical
polynomials $G\in"\R[x_1,...,x_n]\ "$ representable in the form
\begin{equation}G(x)="\sum_{i\in J}h_i(x)F_i(x)\ "=\max_{i\in
J}(h_i(x)+F_i(x))\ ,\label{e4}\end{equation} where $J$ is a finite
non-empty subset of $\{1,2,...,k\}$, and $h_i\in"\R[x_1,...,x_n]\
"$, $i\in J$. The radical $\sqrt{I}$ of an ideal $I$ is the set of
polynomials $G$ such that $"G^m\ "=mG\in I$ for some natural $m$.
We ask the following question:

{\it Given a tropical polynomial $F\in"\R[x_1,...,x_n]\ "$, under
what conditions does $F$ belong to $\sqrt{I(F_1,...,F_k)}$ ?}

The answer, presented in the next section, is given in terms of
"zero loci" of tropical polynomials $F,F_1,...,F_k$, like in the
classical case.

\section{A tropical polynomial Nullstellensatz}

{\it The amoeba} $A(F)\subset\R^n$ of a tropical polynomial $F$ in
$n$ variables is the corner locus of the graph of $F(x)$. This is
a finite $(n-1)$-dimensional polyhedral complex, whose complement
consists of open convex polyhedra. It plays the role of the zero
locus of a tropical polynomial.

\begin{theorem}\label{t1} Let $F,F_1,...,F_k$, $k\ge 1$, be tropical polynomials in $n$
variables. Then $F\in\sqrt{I(F_1,...,F_k)}$ if and only if, for
any connected component $D$ of $\R^n\backslash A(F)$, there is
$1\le i\le k$ such that $D\cap A(F_i)=\emptyset$, and, for each
$j=1,...,n$,
\begin{equation}\frac{\partial F}{\partial x_j}\Big|_D>0\quad\text{as far
as}\quad\frac{\partial F_i}{\partial x_j}\Big|_D>0\
.\label{e1}\end{equation}
\end{theorem}

\begin{proof} {\it 1. Necessity}. The function $mF\big|_D$ is linear.
Hence it must coincide with one of the terms $(h_i+F_i)\big|_D$ in
the expression $\max_i(h_i+F_i)\big|_D$, since otherwise the graph
of the latter function would have a break inside $D$. Next, if
\mbox{$mF\big|_D=(h_i+F_i)\big|_D$}, then both $h_i\big|_D$ and
$F_i\big|_D$ must be linear in view of
\mbox{$A(h_i+F_i)=A(h_i)\cup A(F_i)$}. Thus, $D\cap
A(F_i)=\emptyset$. Observing that
$$m\frac{\partial F}{\partial x_j}\Big|_D=\frac{\partial
h_i}{\partial x_j}\Big|_D+\frac{\partial F_i}{\partial
x_j}\Big|_D\ge\frac{\partial F_i}{\partial x_j}\Big|_D\ ,$$ we
obtain (\ref{e1}).

\medskip

{\it 2. Sufficiency}. Distribute the connected components of
$\R^n\backslash A(F)$ into disjoint subsets $\Pi_i$, $i\in J$,
where $J\subset\{1,...,k\}$, such that, for any $i\in J$ and
$D\in\Pi_i$, we have $A(F_i)\cap D=\emptyset$ and relation
(\ref{e1}).

Fix some $i\in J$. Condition (\ref{e1}) yields that there is $m_1$
such that, for any $m\ge m_1$, one has
$$m\frac{\partial F}{\partial x_j}\Big|_D\ge\frac{\partial
F_i}{\partial x_j}\Big|_D,\quad D\in\Pi_i,\ j=1,...,n\ ,$$ and
hence the gradients of the linear functions
\begin{equation}L_{D,m}:\R^n\to\R,\quad
L_{D,m}\big|_D=mF\big|_D-F_i\big|_D,\quad D\in\Pi_i,\ m\ge m_1\
,\label{e3}\end{equation} have non-negative integral coordinates.

We claim that there is $m_2$ such that, for any $D\in\Pi_i$, in
the complement of the closure of $D$, we have
\begin{equation}mF> L_{D,m}+F_i\quad\text{as far as}\quad m\ge m_1\ .\label{e2}\end{equation} Indeed,
write $F=L+\Phi$, $F_i=L'+\Phi'$, where $L$, $L'$ are linear
functions, and $\Phi$, $\Phi'$ are convex piece-wise linear
functions, vanishing along $D$. Then $L_{D,m}=mL-L'$, and thus,
$mF-L_{D,m}-F_i=m\Phi-\Phi'$. Since $\Phi>0$ and
$|\partial\Phi/\partial x_j|\ge 1$, $j=1,...,n$, outside
$\overline D$, we obtain (\ref{e2}) when $m_2$ exceeds all the
absolute values of the partial derivatives of $\Phi'$.

Define $h_i=\max_{D\in\Pi_i}L_{D,m}$. This is a tropical
polynomial as $m\ge m_1$, and due to (\ref{e3}), (\ref{e2}) it
satisfies
$$(h_i+F_i)\big|_D=mF\big|_D,\quad
(h_i+F_i)\big|_{\R^n\backslash\overline D}\ <\
F\big|_{\R^n\backslash\overline D},\quad D\in\Pi_i,\quad m\ge m_2\
.$$ That is $G="F^m\ "=mF$ satisfies (\ref{e4}) for all
sufficiently large $m$. \end{proof}

\begin{remark}
From the above proof, one can extract an explicit upper bound to
the minimal value of $m$ such that $"F^m\ "\in I(F_1,...,F_k)$.
\end{remark}

\begin{example}
In the case $k=1$, the criterion of Theorem \ref{t1} for
$F\in\sqrt{I(F_1)}$ can be written as \begin{itemize}\item
$A(F)\supset A(F_1)$, and \item for each $j=1,...,n$, $\partial
F/\partial x_j(x)>0$ as far as $\partial F_1/\partial x_j(x)>0$,
$x\not\in A(F)$.\end{itemize} We shall comment on the first
condition. One can assign integer positive weights to the
$(n-1)$-cells of an amoeba so that it will satisfy an equilibrium
condition (see \cite{Mi2}, section 2.1). Taking $m$ greater than
the maximal ratio of weights in $A(F)$ and in $A(F_1)$, we then
subtract the weight of an $(n-1)$-cell $D$ of $A(F_1)$ from the
multiplied by $m$ weights of those $(n-1)$-cells of $A(F)$, whose
interior intersects with $D$. The equilibrium condition persists
after the subtraction due to its linearity, and we again obtain a
balanced complex supported at $A(F)$. By \cite{Mi2}, Proposition
2.4, it defines a convex piece-wise linear function $h_1$, which
provides the relation $mF=h_1+F_1$.
\end{example}

\section{Modifications}

\subsection{A tropical Laurent polynomial Nullstellensatz} A
tropical Laurent polynomial in $n$ variables is a function given
by (\ref{e5}), where $\Omega$ is any non-empty finite subset of
$\Z^n$. Denote the space of tropical Laurent polynomials by
$"L[x_1,...,x_n]\ "$. Correspondingly we define {\it the tropical
Laurent polynomial ideal}
$I^\lau(F_1,...,F_k)\subset"L[x_1,...,x_n]\ "$, generated by
Laurent polynomials $F_1,...,F_k$, as the set of tropical Laurent
polynomials representable in the form (\ref{e4}) with any
non-empty finite $J\subset\{1,...,k\}$, and
$h_i\in"L[x_1,...,x_n]\ "$, $i\in J$.

\begin{theorem}\label{t2}
Let $F,F_1,...,F_k$, $k\ge 1$, be tropical Laurent polynomials in
$n$ variables. Then $F\in\sqrt{I^\lau(F_1,...,F_k)}$ if and only
if, for any connected component $D$ of $\R^n\backslash A(F)$,
there is $1\le i\le k$ such that $D\cap A(F_i)=\emptyset$.
\end{theorem}

\begin{proof} The necessity is established in the same way as in the proof
of Theorem \ref{t1}. The sufficiency follows from Theorem
\ref{t1}, if one "multiplies" $F,F_1,....,F_k$ by suitable
tropical monomials, turning $F,F_1,...,F_k$ into tropical
polynomials, and making all the partial derivatives of $F$
positive. \end{proof}

\subsection{Restricted ideals and restricted Nullstellensatz} In
the space $"\R[x_1,...,x_n]\ "$ introduce {\it the restricted
ideal, generated by tropical polynomials} $F_1,...,F_k$, as
$$I^r(F_1,...,F_k)=\big\{G\in"\R[x_1,...,x_n]\ "\ :$$ $$
G="\sum_{i=1}^kh_iF_i\ "=\max_{1\le i\le k}(h_i+F_i),\quad
h_1,...,h_k\in"\R[x_1,...,x_n]\ "\big\}\ .$$ Notice that
$I^r(F_1,...,F_k)\subset I(F_1,...,F_k)$, but they may differ, for
example, \mbox{$x\not\in I^r(x,"x+1\ ")$}.

\begin{theorem}\label{t3}
Let $F,F_1,...,F_k$, $k\ge 1$, be tropical polynomials in $n$
variables. Then $F\in\sqrt{I^r(F_1,...,F_k)}$ if and only if the
following conditions hold:
\begin{enumerate}\item[(i)] for any connected component $D$ of $\R^n\backslash A(F)$,
there is $1\le i\le k$ such that $D\cap A(F_i)=\emptyset$ and
(\ref{e1}) is fulfilled; \item[(ii)] there is $m_0$ such that
$m\Del(F)$ contains a translate of each Newton polytope
$\Del(F_1),...,\Del(F_k)$ as far as $m\ge m_0$.
\end{enumerate}
\end{theorem}

\begin{remark}
Notice that condition (ii) always holds when $\dim\Del(F)=n$.
\end{remark}

\begin{proof} {\it 1. Auxiliary statement}. Let $G,H\in"\R[x_1,...,x_n]\ "$.
We claim that
\begin{itemize}\item if $G\ge H$, then $\Del(G)\supset\Del(H)$;
\item if $\Del(G)\supset\Del(H)$, then there is $c\in\R$ such that
$G\ge H+c$.\end{itemize}

Represent $\Del(G)$ as intersection of finitely many halfspaces.
Pick one of these halfspaces and apply an integral-affine
automorphism $Q$ of $\R^n$, taking the halfspace to $x_n\ge0$.
Since $G\circ Q^{-1}\ge H\circ Q^{-1}$, the Newton polytope
$\Del(H\circ Q^{-1})=Q(\Del(H))$ cannot contain points with
negative $n$-th coordinate. Indeed, otherwise we would have that,
for $x_1,...,x_{n-1}=\const$, $x_n\to-\infty$, the function
$G(Q^{-1}(x))$ does not increase, whereas $H(Q^{-1}(x))$ tends to
$+\infty$. Running over all halfspaces forming $\Del(G)$, we
conclude that $\Del(G)\supset\Del(H)$.

For the second statement, write
$$G(x)=\max_{\omega\in\Del(G)\cap\Z^n}(\langle
x,\omega\rangle+a_\omega),\quad
H(x)=\max_{\omega\in\Del(H)\cap\Z^n}(\langle
x,\omega\rangle+b_\omega)\ .$$ Then one can take
$$c=\min_{\omega\in\Del(H)\cap\Z^n}(a_\omega-b_\omega)\ .$$

{\it 2. Necessity}. We have to prove only (ii). Observing that
$\Del(h_i+F_i)$ is the convex hull of few translates of
$\Del(F_i)$, we derive (ii) from the above auxiliary
statement.

{\it 3. Sufficiency}. A translate of $\Del(F_i)$ is the Newton
polytope of a tropical polynomial $"h_iF_i\ "$, where $h_i$ is a
tropical monomial (i.e, a linear function). Then, according to the
auxiliary statement, $mF_i\ge (h_i+c_i)+F_i$ for sufficiently
large $m$ and certain constant $c_i$. Then, using Theorem
\ref{t1}, we obtain $$mF=\max_{i\in
J}(h_i+F_i)=\max\left\{\max_{i\in J}(h_i+F_i),\ \max_{i\not\in
J}(h_i+c_i+F_i)\right\}\ .$$\end{proof}

\subsection{Nullstellensatz for convex piece-wise linear
functions of finite type} A function given by (\ref{e5}), where
$\Omega$ is a finite subset of $\R^n$ we call {a convex piece-wise
linear function of finite type} in $n$ variables. Denote the space
of such functions by $"PL(x_1,...,x_n)\ "$. As in the polynomial
case, we define amoebas and finitely generated ideals
$I^\fun(F_1,...,F_k)$ in $"PL(x_1,...,x_n)\ "$, and obtain a
corresponding Nullstellensatz:

\begin{theorem}\label{t4}
Let $F,F_1,...,F_k\in"PL(x_1,...,x_n)\ "$. Then
$F\in\sqrt{I^\fun(F_1,...,F_k)}$ if and only if, for any connected
component $D$ of $\R^n\backslash A(F)$, there is $1\le i\le k$
such that $D\cap A(F_i)=\emptyset$.
\end{theorem}

Proof coincides with that of Theorem \ref{t2}.

\subsection{Nullstellensatz for polynomials over an extended
tropical semiring} In \cite{Iz}, the second author introduced an
extension $\T$ of the tropical semiring $\R$, obtained by adding
one more copy of $\R$, which we denote by $\R^\nu$ and its
elements by $a^\nu\in\R^\nu$ as $a\in\R$, and equipped with the
tropical operations
$$"a+b\ "=\begin{cases}\max\{a,b\},\ &\text{if}\ a\ne b,\\
a^\nu,\ &\text{if}\ a=b,\end{cases}\qquad"a^\nu+b^\nu\
"=(\max\{a,b\})^\nu\ ,$$
$$"a+b^\nu\ "=\begin{cases}a,\ &\text{if}\ a>b,\\ b^\nu,\ &\text{if}\
a\le b,\end{cases}$$ $$"ab\ "=a+b,\quad"a^\nu b^\nu\
"=(a+b)^\nu,\quad "ab^\nu\ "=(a+b)^\nu$$ for all $a,b\in\R$. There
is a natural epimorphism of tropical semirings:
$$\pi:\T\to\R,\quad \pi(a)=a,\ \pi(a^\nu)=a,\quad\text{for all}\
a\in\R\ ,$$ which induces epimorphisms $\pi_*$ of the tropical
polynomial and tropical Laurent polynomial rings. The above
variants of Nullstellensatz can easily be translated to the case
of base semiring $\T$. We present here such a translation of
Theorem \ref{t1}.

For a polynomial $F="\sum_{\omega\in\Omega}c_\omega x^\omega\
"\in"\T[x_1,...,x_n]\ "$, where $\Omega$ is a non-empty finite set
of points with non-negative integral coordinates, we introduce an
amoeba $A(F):=A(\pi_*F)\subset\R^n$. The restriction of $\pi_*F$
to a connected component $D$ of $\R^n\backslash A(F)$ is a linear
function $\langle\omega,x\rangle+\pi(c_\omega)$ for some
$\omega=\omega(D)\in\Omega$. So we divide the set of the connected
components of $\R^n\backslash A(F)$ into two disjoint subsets
$\Pi(F)$ and $\Pi^\nu(F)$, letting $D\in\Pi(F)$ or
$D\in\Pi^\nu(F)$ according as $c_{\omega(D)}\in\R$ or
$c_{\omega(D)}\in\R^\nu$.

\begin{theorem}\label{t5}
Let $F,F_1,...,F_k\in"\T[x_1,...,x_n]\ "$, $k\ge 1$. Then
$F\in\sqrt{I(F_1,...,F_k)}$ if and only if
\begin{enumerate}\item[(i)] for any $D\in\Pi(F)$, there is $1\le i\le k$ such that $D$ is contained in some component
$D_i\in\Pi(F_i)$, and, for each $j=1,...,n$, relation (\ref{e1})
holds true; \item[(ii)] for any $D\in\Pi^\nu(F)$, there is $1\le
i\le k$ such that $D\cap A(F_i)=\emptyset$, and, for each
$j=1,...,n$, relation (\ref{e1}) holds true.
\end{enumerate}
\end{theorem}

\begin{proof} The necessity part is immediate in view of the proof
of Theorem \ref{t1} and the fact that $"\R\R^\nu\ "=\R^\nu$. To
prove the sufficiency, we construct coefficients $h_i$ of the
expansion (\ref{e4}) to be in $"\R[x_1,...,x_n]\
"\subset"\T[x_1,...,x_n]\ "$ by the following modification of the
procedure from the proof of Theorem \ref{t1}:
\begin{itemize}\item if $D\in\Pi(F)$, and $1\le i\le k$ satisfies
condition (i), then we put $h_i\big|_D=mF\big|_D-F_i\big|_xD$;
\item if $D\in\Pi^\nu(F)$, and $1\le i\le k$ satisfies condition
(ii), then we put
$h_i\big|_D=m\cdot\pi_*F\big|_D-\pi_*F_i\big|_D$. \end{itemize}
Choosing a sufficiently large $m$, we complete the proof as that
of Theorem \ref{t1}. \end{proof}

\begin{corollary}
Let $F,F_1,...,F_k\in"\T[x_1,...,x_n]\ "$, $k\ge 1$, and
$I(F_1,...,F_k)\supset"\R^\nu[x_1,...,x_n]\ "$. Then
$F\in\sqrt{I(F_1,...,F_k)}$ if and only if condition (i) of
Theorem \ref{t5} holds true.
\end{corollary}

Indeed, the relation $I(F_1,...,F_k)\supset"\R^\nu[x_1,...,x_n]\
"$ means that the set of generators $\{F_1,...,F_k\}$ contains a
constant polynomial $F_i=a^\nu$, $a\in\R$, and thus, condition
(ii) always holds.

\bibliographystyle{amsplain}

\end{document}